\documentclass[12pt,twoside]{article}

\usepackage[T1]{fontenc}
\usepackage[cp1250]{inputenc}

\usepackage[arrow, matrix, curve]{xy}
\usepackage{amsmath,amsthm,amssymb}
\usepackage{amscd,graphics}
\usepackage[dvips]{graphicx}
\usepackage{makeidx}

\newenvironment{Proof}{\textbf{Proof.}}{$\qquad \blacksquare$\par}
\newenvironment{Proof of}[1]{\textbf{Proof #1.}}{$\qquad \blacksquare$\par}
\newenvironment{Item}[1]{\par #1 }{\par}

\newcommand{\Ker}{\textrm{Ker}\,}
\newcommand{\hull}{\textrm{hull}\,}

\newcommand{\B}{\mathcal B}
\newcommand{\M}{\mathcal M}

\newcommand{\NN}{\mathcal N}

\newcommand{\al}{\alpha}

\newcommand{\A}{\mathcal A}

\newcommand{\Z}{\mathbb Z}
\newcommand{\N}{\mathbb N}

\newtheorem{thm}{Theorem}[section]

\newtheorem{prop}[thm]{Proposition}
\newtheorem{cor}[thm]{Corollary}
\theoremstyle{definition}

\newtheorem{defn}[thm]{Definition}
\newtheorem{ex}[thm]{Example}
\newtheorem{rem}[thm]{Remark}

 \makeindex

\begin{document}

 \thispagestyle{empty}

 \begin{center}
{\bfseries \large \textsc{Crossed product by an arbitrary endomorphism}}

\bigskip
B. K.  Kwa\'sniewski, \ \  A. V. Lebedev

\end{center}

\begin{abstract} Starting
from an arbitrary endomorphism $\delta$ of a unital $C^*$-algebra
$\A$ we construct a crossed product. It is shown that the natural construction depends not only on the $C^*$-dynamical system $(\A,\delta )$ but also on the choice of an ideal $J$ orthogonal to $\Ker \delta $.
\end{abstract}

\medbreak


\medbreak
\vspace{5mm}

\tableofcontents

\section*{Introduction}

The crossed product of a $C^*$-algebra $\A$ by an automorphism
$\delta:\A\to \A$ is defined as a universal $C^*$-algebra
generated by a copy of $\A$ and a unitary element $U$ satisfying
the relations
$$
\delta(a)= Ua U^*,\qquad \delta^{-1}(a)=U^*aU, \ \ \ \ a\in \A.
$$
On  one hand, algebras arising in this way (or their versions adapted to actions
of groups of automorphisms) are very well understood and  became
a part of a $C^*$-folklore \cite{Kadison}, \cite{Pedersen}. On the other hand,
it is  very symptomatic  that, even though  the first attempts on  generalizing
this kind  of constructions to endomorphisms  go  back to 1970s,  articles introducing
different  definitions of the related object appear almost  continuously until the present-day, see,
for example, \cite{CK}, \cite{Paschke},
\cite{Stacey}, \cite{Murphy},  \cite{exel1}, \cite{exel2},
\cite{kwa}, \cite{Ant-Bakht-Leb}.
This phenomenon is caused by   very fundamental problems one has to face when dealing with crossed
products by endomorphisms. Namely, one has to answer the following questions:

\begin{itemize}
    \item[(i)] What relations should the element  $U$  satisfy?
        \item[(ii)] What should be used in place of $\delta^{-1}$?
\end{itemize}

It is important that in spite of the  substantial freedom of choice (in answering the foregoing questions),
all the above listed papers  do however have a certain nontrivial intersection. They  mostly agree, and simultaneously
boast their greatest successes, in the  case when dynamics is
implemented by  monomorphisms with  a hereditary range. In view of
the recent articles \cite{Bakht-Leb}, \cite{Ant-Bakht-Leb}, \cite{kwa3}, this coincidence is completely
understood. It is shown in \cite{Bakht-Leb} that in the case of
monomorphism with hereditary range there exists a unique
non-degenerate transfer operator $\delta_*$ for $(\A,\delta)$,
called by authors of  \cite{Bakht-Leb} a \emph{complete transfer
operator}, and the theory goes smooth with  $\delta_*$ as it takes
over the role classically played by $\delta^{-1}$. The $C^*$-dynamical
systems of this sort will be  called {\em partially reversible}.
\par
If the pair $(\A,\delta)$ is of the above described type
then $\A$ is called a {\em coefficient algebra}.
This notion was introduced in \cite{Leb-Odz} where its investigation and
relation to the extensions of $C^*$-algebras by partial isometries was clarified.
Further in
 \cite{Bakht-Leb} a certain criterion for a
$C^*$-algebra to be a coefficient algebra associated with a given
endomorphism was obtained. On the base of the results of these papers
one naturally arrives at the construction of a certain crossed
product which was implemented in \cite{Ant-Bakht-Leb}.
It was also observed in \cite{Ant-Bakht-Leb} that in the most
natural situations the coefficient algebras arise as a result of
a certain extension procedure on the initial $C^*$-algebra.
Since the crossed product is (should be) an extension of the
initial $C^*$-algebra one can consider the construction of an
appropriate coefficient algebra as one of the most important
intermediate steps in the procedure of construction of the crossed product itself.
Thus one arrives at the next natural problem: can we extend a
 $C^*$-dynamical system associated with an arbitrary endomorphism
to a {\em partially reversible}  $C^*$-dynamical system? In the
{\em commutative} $C^*$-algebra situation the corresponding
procedure and the explicit description of maximal ideals of the
arising $C^*$-algebra is given in \cite{maxid}. On the base of
this construction the general construction of the crossed product
associated to an arbitrary endomorphism of a {\em commutative}
$C^*$-algebra is presented in \cite{kwa}. Further the general
construction of an extension of a $C^*$-dynamical system
associated with an arbitrary endomorphism to a  partially
reversible $C^*$-dynamical system   is worked out in
\cite{kwa3}, \cite{kwa4}. Therefore the mentioned results of
\cite{Bakht-Leb}, \cite{Ant-Bakht-Leb}, \cite{kwa3}, \cite{kwa4}, give us the key to
construct a general crossed product starting from a
$C^*$-dynamical system associated with an arbitrary endomorphism,
and this is the theme of the present article,  see also Remark  \ref{remark something}.
\par
The main novelty here is
\begin{itemize}
    \item[(i)] the explicit description of the crossed product based on the
worked out matrix calculus presented in Section \ref{2},
\end{itemize}
and an observation of (in a way unexpected) phenomena  that
\begin{itemize}
\item[(ii)] in general the universal construction of the crossed product
depends not only on the algebra $\A$ and an endomorphism $\delta$ one starts with but
also on the choice of an ({\em arbitrary}) singled out ideal $J$
orthogonal to the kernel of $\delta$ (see Section \ref{1}).
\end{itemize}
So in fact we have a {\em variety} of crossed products depending
on $J$.
\par
The paper is organized as follows.
\par
In the first section we recall a notion of  covariant representation  of $C^*$-dynamical system
which appear  in a similar or identical form, for instance, in  \cite{kwa4}, \cite{Leb-Odz}, \cite{kwa},
   \cite{exel2}, \cite{Murphy}, \cite{Stacey}. We split  a class of covariant representation
   into subclasses according to a certain ideal  they determine, and explain a role
   this distinction plays in the present article.
 Section 2 contains  a matrix calculus  which  describes an algebraic structure
 of the crossed product and  in Section 3  we calculate the appropriate norms.
 Gathering the facts from the previous sections we define  crossed product in Section 4.
 In the final section we include theorems concerning   representations of the newly defined object.

\section{Covariant representations and orthogonal ideals}
\label{1}

In this section we discuss interrelations between covariant representations
and orthogonal ideals. In particular,  we indicate (in Remark  \ref{remark something})
a definition of the crossed product.

Let $(\A,\delta)$ be a   pair consisting of a $C^*$-algebra $\A$, containing an identity and
an endomorphism $\delta : \A \to \A$. Throughout the paper the pair $(\A,\delta)$ will be called
a {\em $C^*$-dynamical system}.

\begin{defn}\label{kowariant  rep defn}
Let $(\A,\delta)$ be a $C^*$-dynamical system. A triple
$(\pi,U,H)$ consisting of a non-degenerate
faithful representation $\pi:\A\to \B(H)$ on a Hilbert space $H$
and a partial isometry $U\in \B(H)$ satisfying
$$
U\pi(a)U^* =\pi(\delta(a)),\qquad a \in \A, \qquad U^*U \in \pi(\A)',
$$
 will be called
a \emph{covariant representation} of  $(\A,\delta)$.
\end{defn}

The first named author  described in \cite{kwa4} a procedure of extending any $C^*$-dynamical system
$(\A,\delta)$   up to a system $(\A^+,\delta^+)$ with a property that  a kernel of $\delta^+$ is unital.
Moreover, the resulting system $(\A^+,\delta^+)$ is in a sense the smallest  extension of $(\A,\delta)$
possessing that property, see \cite{kwa4}.
Let us now slightly generalize this construction, which will be essential for  our future purposes.
\\
Let $(\A,\delta)$  be a $C^*$-dynamical system and denote by $I$ the kernel of $\delta:\A\to \A$.
We shall say that an {\em ideal} $J$ in $\A$ is  \emph{ orthogonal to} $I$ if
$$
I\cap J =\{0\}.
$$
There exists the biggest ideal orthogonal  to $I$
(in the sense that it contains all other orthogonal ideals to $I$)
and we shall denote it by $I^\bot$. This ideal  could be defined explicitly. Namely if we let
$$
\hull(I)=\{x\in {\rm Prim} A, \ x\supset I \}
$$
then we have
$$
I= \bigcap_{x\in \hull(I)}\, x, \ \qquad  \ \ I^\bot = \bigcap_{x\in {\rm Prim} A\setminus \hull(I)}\,  x.
$$
Alternative explicit definition of $I^\bot$ can be found  in  \cite{kwa4}.\\
 Our construction will depend on the  choice of an ideal
orthogonal to $I$. Let us single out  a certain ideal $J$ which is
orthogonal to $I$ in $\A$, that is
$$
\{0\} \,\subset \,J\, \subset\,  I^\bot.
$$
By    $\A_J$  we denote the direct sum of quotient algebras
$$
\A_J=\big(\A/I\big) \oplus \big(\A/J\big),
$$
and  we set $\delta_J:\A_J\to \A_J$ by the formula
\begin{equation}
\label{d_J}
 \A_J\ni \big((a +I)\oplus (b
+J)\big)\stackrel{\delta_J}{\longrightarrow}(\delta(a) +I)\oplus
(\delta(a) +J)\in \A_J.
\end{equation}
Routine verification shows that
\begin{equation}
\label{kerd_J} \Ker \delta_J = (0,\, \A/J )
\end{equation}

 Since $I={\rm Ker}\delta$ it follows that an element
$\delta(a)$ does not depend  on the choice of   a representative
of $a +I$ and so the mapping $\delta_J$ is well defined. Clearly,
$\delta_J$ is an endomorphism and its  kernel is unital with the
unit of the form $(0 + I) \oplus (1 +J)$.
\\
Moreover, the $C^*$-algebra $\A$ embeds  into $C^*$-algebra $\A_J$ via
\begin{equation}
\label{A_in} \A \ni a \longmapsto \big(a +I\big)\oplus \big(a +
J\big)\in \A_J.
\end{equation}

 Since $I\cap J = \{0\}$ this mapping is injective  and we shall
treat $\A$ as the corresponding subalgebra of $\A_J$. Under this
identification $\delta_J$ is an extension of $\delta$.
\\
The main motivation of the preceding construction is the following statement.
\begin{prop}\label{motivation prop}
Let $(\pi,U,H)$ be a covariant representation of $(\A,\delta)$ and let
$$
I=\{ a\in \A: (1-U^*U) \pi(a)=\pi(a)\}, \qquad J=\{ a\in \A: U^*U \pi(a)=\pi(a)\}.
$$
Then  $I$ is the kernel of $\delta$ and $J$ is an ideal orthogonal to $I$. Moreover,
if $(\A_J,\delta_J)$ is the extension of $(\A,\delta)$ constructed above,
then  $\pi$ uniquely extends to the isomorphism $\widetilde{\pi}:\A_J\to C^*(\pi(\A),U^*U)$ such that
$(\widetilde{\pi},U,H)$ is a covariant representation of $(\A_J,\delta_J)$.
Namely $\widetilde{\pi}$ is given by
\begin{equation}\label{extension of pi to J-algebra}
\widetilde{\pi}(a + I \oplus b + J)=  U^*U \pi(a) + (1- U^*U) \pi(b), \qquad a,b \in \A.
\end{equation}
\end{prop}
\begin{Proof}
Observe that
$$
U (U^*U \pi(a))U^* = U  \pi(a)U^* =U (U^*U \pi(a))U^*
$$
and
$$
U^* (U \pi(a)U^*)U =U^*U \pi(a).
$$
Therefore the mappings
$$
U(\cdot )U^*:U^*U \pi(a)\mapsto\pi(\delta(a)),
$$
$$
U^*(\cdot )U:\pi(\delta(a))\mapsto U^*U \pi(a)
$$
establish an isomorphism
 $U^*U \pi(\A)\cong \pi(\delta(\A)), \ (U^*U\pi(a) \mapsto
 \pi(\delta(a))).$ Thus
$$
\delta(a)=0  \Longleftrightarrow \pi(\delta(a))=0
\Longleftrightarrow U^*U\pi(a)=0 \Longleftrightarrow
(1-U^*U)\pi(a)=\pi(a).
$$
Which means that   $I$ is the kernel of $\delta$.

It is clear that $J$ is an ideal orthogonal to $I$.
If $\widetilde{\pi}:\A_J\to C^*(\pi(\A),U^*U)$ is onto and
$(\widetilde{\pi},U,H)$ is a covariant representation of
$(\A_J,\delta_J)$ then by \cite[Proposition 1.5]{kwa4}   we have
$$
U^*U = \widetilde{\pi}(1 + I \oplus 0 + J).
$$
Thus $\widetilde{\pi}$ is of the form \eqref{extension of pi to
J-algebra}. Conversely, formula \eqref{extension of pi to
J-algebra} defines an isomorphism which follows from the
definitions of $I$ and $J$. It is readily checked  that
$(\widetilde{\pi},U,H)$ is a covariant representation of
$(\A_J,\delta_J)$.
\end{Proof}
\begin{rem} The larger ideal  $J$ is the smaller $\A_J$ is.
Namely, since $  {\rm Ker}{\delta_J} = (0,\, \A/J ) $ (recall
(\ref{kerd_J})) it follows that
 as $J$ varies from $\{0\}$ to
$I^\bot$ the kernel of $\delta_J$  varies from $(0,\, \A )$ to
$(0,\,\A/I^{\bot})$. On the other hand, in view of (\ref{d_J}) and
(\ref{A_in}),
 the image of
$\delta_J$ is always isomorphic to $\A/I$ (the first  summand of
$\A_J$) and thereby  it  does not depend  on the choice of $J$.
The case when $J=I^\bot$ was considered in \cite{kwa4}, see
discussion below.
\end{rem}

Proposition \ref{motivation prop} makes it natural to introduce the following definition.

\begin{defn}\label{kowariant  rep defn 2} Let $(\pi,U,H)$ be a covariant representation
$(\pi,U,H)$ of $(\A,\delta)$ and $J$  an ideal in $\A$ orthogonal to kernel of $\delta$.
If   $(\pi,U,H)$  and  $J$ are combined with each other via equality
$$
J=\{ a\in \A: U^*U \pi(a)=\pi(a)\}
$$
then we shall say that
 $(\pi,U,H)$  is a  \emph{covariant representation associated with the ideal } $J$
 and also that $J$ is the \emph{ideal associated to covariant representation} $(\pi,U,H)$.
\end{defn}
\begin{rem}\label{remark something}
We stress that in view of Proposition  \ref{motivation prop} and results of \cite{kwa4}
and \cite{Ant-Bakht-Leb} one could  proceed to  defining   crossed product   associated with an ideal $J$,
right away. Namely, if we fix an ideal $J$ orthogonal to kernel of $\delta$ then the algebra
$$
C^*(\pi(\A),U^*U) \cong \A_J
$$
does not depend on the choice of a covariant representation $(\pi,U,H)$ (as long as it is associated with $J$).
Furthermore, it was shown in \cite{kwa4} that the algebra
$$
\B=C^*(\bigcup_{n\in \N}U^{*n}\pi(\A)U^{n})
$$
can be  described in terms of the system $(\A_J,\delta_J)$, and  therefore it also
does not depend on the choice of $(\pi,U,H)$.  Obviously we have
$$
U\B U^*\subset \B, \qquad U^*\B U\subset \B, \qquad U^*U\in Z(\B)
$$
which means that $\B$ is a \emph{coefficient algebra} for $C^*(\pi(\A), U)$ and
the $C^*$-dynamical system on $\B$, with dynamics implemented by $U$, is \emph{partially reversible}.
Thus  we can  apply  the crossed product elaborated in  \cite{Ant-Bakht-Leb} to the extended
partially reversible system on $\B$ (which is uniquely determined by $(\A,\delta)$ and $J$).
\\
Summing up, the construction of crossed product involve essentially two steps: extending an
\emph{irreversible} system on $\A$  up to a \emph{partially reversible} system on $\B$,
and then attaching the crossed product to the extended system. Such an approach was explored
in the commutative case in \cite{kwa}.
In contrast to this two (or even three) step construction,   in the present paper, however,
we propose a direct approach which provides us with a one more interesting, integral point of view.
\end{rem}

In view of Proposition \ref{motivation prop} and \cite[Theorem 3.7]{kwa4}, cf. above remark,
we have the following highly non-trivial statement.
\begin{thm}\label{existance theorem}
Let $J$ be   an arbitrary ideal in $\A$ having  zero intersection with the kernel of $\delta$.
There exists a   covariant representation associated with the ideal $J$.
\end{thm}

The first named author introduced in \cite{kwa4} a notion of a strict covariant representation
where by a \emph{strict covariant representation of }$(\A,\delta)$ he meant
a covariant representation  $(\pi,U,H)$   such that
\begin{equation}\label{stricct covariant eqation}
U^*U=1-s\lim_{\lambda_\in \Lambda} \pi(\mu_\lambda)
\end{equation}
where $\{\mu_\lambda\}_{\lambda\in \Lambda}$ is an approximate unit in the kernel $I$ of $\delta$.
It follows that if  $(\pi,U,H)$ satisfies \eqref{stricct covariant eqation}
it is associated with the ideal $I^\bot$, cf.  \cite[Remark 1.11]{kwa4}. However the converse statement is not true.
\begin{ex}
Let $\A=c$ be the algebra of all converging sequences (numbered by $1,2,3,4,...$)
and let $\delta$ be an endomorphism of $\A$ which maps a sequence $(a(1),a(2),...)$
onto a constant one $(a(\infty), a(\infty),...)$ where $a(\infty)=\lim_{k\to \infty}a(k)$.
 Then the kernel $I=c_0$ of $\delta$ is the algebra of all sequences converging to zero.
 Let $H=l_2(\Z)$ and  define  $\pi:\A\to L(H)$  by
$$
(\pi(a)h)(k)=
\begin{cases}
 a(k)h(k), & k > 0\\
 a(\infty)h(k), & k\leq 0
\end{cases}.
$$
Let us consider  the family of co-isometries $U_n\in L(H)$,  $n=0,1,2,...$, where
$$
(U^*_nh)(k)=
\begin{cases}
 h(-2k+1-n), & k > 0\\
 h(2k-n), & k\leq 0
\end{cases}.
$$
Then for every $n$ the triple $(\pi,U_n,H)$ is a covariant representations associated with $I^\bot$.
However,  $(\pi,U_n,H)$ satisfies \eqref{stricct covariant eqation} if and only if $n=0$.
\end{ex}
\section{Matrix calculus for covariance algebras}
\label{2}

In this section we describe a certain matrix calculus that will be an algebraic framework for
constructing our principal object --- the crossed product.

 We denote by $\M(\A)$ the set of infinite matrices
$\{a_{ij}\}_{i,j\in\N} $  with entries in $\A$ indexed by pairs of
natural numbers such that the coefficients $a_{ij}$ are of the
form
$$ a_{ij}\in \delta^i(1)\A \delta^j(1), \qquad i, j\in \N,$$
and there is at most finite number of $a_{ij}$ which are non-zero.
\\
We will take advantage of this standard matrix notation when
defining operations on $\M(\A)$ and investigating a natural homomorphism from $\M(\A)$ to the
covariance algebras generated by covariant representations. However,  when investigating
norm of elements in covariance algebras, it more handy to index the  entries of an element
in $\M(\A)$ by a pair consisting of a natural number and an integer.
Hence  we shall paralellely use two notations concerning matrices in $\M(\A)$.
Namely, we presume the following identifications
 $$
a_{ij}=a^{(j-i)}_{\min\{i,j\}},\quad i,j \in \N,\quad  \qquad a_{n}^{(k)}=\begin{cases}
a_{n,k+n},& k\geq 0\\
a_{n-k,n},& k < 0
\end{cases}\quad n \in \N,\,\, k\in \Z,
$$
under which we have two equivalent matrix presentations
$$
 \left(
\begin{array}{cc c c}
  a_{00}      &  a_{01}   &  a_{02}  &   \cdots  \\
  a_{10}     &  a_{11}   &  a_{12}  &   \cdots  \\
  a_{20}     &  a_{21}  &  a_{22}  &   \cdots  \\
    \vdots       &   \vdots     &   \vdots    &    \ddots
\end{array}\right)
=
\left(
\begin{array}{cc c c}
  a_0^{(0)}      &  a_0^{(1)}   &  a_0^{(2)}  &   \cdots  \\
  a_0^{(-1)}     &  a_1^{(0)}   &  a_1^{(1)}  &   \cdots  \\
  a_0^{(-2)}     &  a_1^{(-1)}  &  a_2^{(0)}  &   \cdots  \\
    \vdots       &   \vdots     &   \vdots    &    \ddots
\end{array}\right).
$$
The use of one of this   conventions  will always be clear from the context.
\\
 We define the
addition, multiplication by scalar, and
involution on $\M(\A)$  in a natural manner. Namely, let $a=\{a_{ij}\}_{i,j\in\N}$ and $b=\{b_{ij}\}_{i,j\in\N}$.
We put
\begin{gather}\label{add1}
 (a+b)_{m,n}=a_{m,n}+b_{m,n},\\[6pt]
\label{mulscal1}
 (\lambda a)_{m,n}=\lambda a_{m,n}\\[6pt]
\label{invol1}
 (a^*)_{m,n}=a_{n,m}^*.
\end{gather}
 Moreover,  we introduce a convolution multiplication "$\star$" on
$\M(\A)$, which is a reflection of the operator multiplication in covariance algebras.
 We set
\begin{equation}\label{star1}
a\star b= a\cdot\sum_{j=0}^\infty  \Lambda^j(b)+ \sum_{j=1}^\infty \Lambda^j(a)\cdot b
\end{equation}
where $\cdot$ is the standard multiplication of matrices and  mapping  $\Lambda:
\M(\A)
\rightarrow \M(\A)$ is defined
 to act as follows: $\Lambda(a)_{ij}=\delta(a_{i-1,j-1})$, for $i,j> 1$, and $\Lambda(a)_{ij}=0$ otherwise,
 that is $\Lambda$ assumes the following shape
\begin{equation}\label{Lambda}
\Lambda(a)= \left(
\begin{array}{cc c c c }
 0      &         0           &           0         &            0       & \cdots  \\
 0      &  \delta(a_{00})  &  \delta(a_{01})  &  \delta(a_{02}) & \cdots  \\
 0      &  \delta(a_{10}) &  \delta(a_{11})  &  \delta(a_{12}) & \cdots  \\
 0      &  \delta(a_{20}) &  \delta(a_{21}) &  \delta(a_{22}) & \cdots  \\
\vdots  &       \vdots        &       \vdots        &       \vdots       & \ddots
\end{array}\right).
\end{equation}
\begin{prop}
The set $\M(\A)$ with operations \eqref{add1},
\eqref{mulscal1}, \eqref{invol1}, \eqref{star1} becomes an algebra with involution.
\end{prop}
\begin{Proof}
The only thing we show is the associativity of the multiplication \eqref{star1}, the rest is straightforward.
For that purpose we note that $\Lambda$ preserves the standard matrix multiplication and thus we have
$$
a\star (b \star c)= a\star  \Big(b\cdot\sum_{j=0}^\infty  \Lambda^j(c)+ \sum_{j=1}^\infty \Lambda^j(b)\cdot c\Big)
$$
$$
=a\sum_{k=0}^\infty  \Lambda^k\Big( b\sum_{j=0}^\infty  \Lambda^j(c)+
\sum_{j=1}^\infty \Lambda^j(b) c \Big)+ \sum_{k=1}^\infty \Lambda^k(a) \Big( b\sum_{j=0}^\infty  \Lambda^j(c)+
\sum_{j=1}^\infty \Lambda^j(b) c
\Big)$$
$$
=\sum_{k,j=0}^\infty  a \Lambda^k( b)\Lambda^{j}(c) +  \sum_{k=1,j=0}^\infty \Lambda^{k}(a) b\cdot \Lambda^j(c) +
\sum_{k,j=1}^\infty \Lambda^{k}(a)\Lambda^{j}(b) c
$$
$$
=\Big(a\sum_{k=0}^\infty  \Lambda^k(b)+ \sum_{k=1}^\infty \Lambda^k(a) b\Big)\sum_{j=0}^\infty \Lambda^j(c) +
\sum_{j=1}^\infty \Lambda^j \Big(a\sum_{k=0}^\infty  \Lambda^k(b)+ \sum_{k=1}^\infty \Lambda^k(a) b\Big) c
$$
$$
  \Big(a\cdot\sum_{k=0}^\infty  \Lambda^k(b)+ \sum_{k=1}^\infty \Lambda^k(a)\cdot b\Big)\star c= (a\star b) \star c.
$$
\end{Proof}
We embed $\A$ into $\M(\A)$ by identifying element $a\in \A$ with
a matrix $\{a_{mn}\}_{m,n\in \N}$ where $a_{00}=a$ and $a_{mn}=0$
if $(m,n)\neq (0,0)$.  We also define a partial isometry
$u=\{u_{mn}\}_{m,n\in \N}$ in $\M(\A)$ such that
$u_{01}=\delta(1)$ and $u_{mn}=0$ if $(m,n)\neq (0,1)$. In other
words, we adopt the following notation
\begin{equation}
\label{notation stupid} u= \left(
\begin{array}{cc c c}
  0      &  \delta(1)   &  0  &   \cdots  \\
  0     &  0   &  0  &   \cdots  \\
  0     & 0  &  0  &   \cdots  \\
    \vdots       &   \vdots     &   \vdots    &    \ddots
\end{array}\right)\quad \mathrm{and \,\,} \quad
a=\left(\begin{array}{cc c c}
  a      &  0     &  0    & \cdots  \\
  0      &  0 &  0    & \cdots  \\
  0     &  0    &  0   & \cdots  \\
     \vdots                 &    \vdots                &      \vdots             &\ddots
\end{array}\right), \quad \mathrm{for \,\,} a\in \A.
\end{equation}
 One can easily check that the $^*$-algebra $\M(\A)$ is generated
by $u$ and $\A$. Furthermore, for every $a\in \A$ we have
$$
 u\star a \star u^* =\delta(a) \quad\textrm{ and } \quad  u^*\star a \star u = \left(
\begin{array}{cc c c}
  0      &  0    &  0  &   \cdots  \\
  0     &  \delta(1)a\delta(1)   &  0  &   \cdots  \\
  0     & 0  &  0  &   \cdots  \\
    \vdots       &   \vdots     &   \vdots    &    \ddots
\end{array}\right).
 $$
 \begin{prop}\label{dense subalgera structure prop}
Let $(\pi,U,H)$ be  a covariant representation of $(\A,\delta)$.
Then there exists a  unique $^*$-homomorphism $\Psi_{(\pi, U)}$
from $\M(\A)$ onto a $^*$-algebra $C_0^*(\pi(\A),U)$ generated by
$\pi(\A)$ and $U$, such that
$$
\Psi_{(\pi, U)}(a)=\pi(a),\quad a\in \A,\qquad \Psi_{(\pi, U)}(u)=U.
$$
Moreover, $\Psi_{(\pi, U)}$ is given by the formula
\begin{equation}\label{Psi form eq}
\Psi_{(\pi, U)}( \{a_{m,n}\}_{m,n\in\N})=\sum_{m,n=0}^\infty U^{*m}\pi(a_{m,n}) U^n,
\end{equation}
and thus $C_0^*(\pi(\A),U)=\{\sum_{m,n=0}^\infty U^{*m}\pi(a_{m,n}) U^n: \{a_{m,n}\}_{m,n\in\N}\in \M(\A)\}$.
\end{prop}
\begin{Proof} It is clear that $\Psi_{(\pi, U)}$ has to satisfy  form  \eqref{Psi form eq}.
Thus it is enough to check that $\Psi_{(\pi, U)}$ is a $^*$-homomorphism, and in fact we only need to show
that $\Psi_{(\pi, U)}$ is multiplicative as the rest is obvious. For that purpose
let us fix two matrices $a=\{a_{m,n}\}_{m,n\in\N},b=\{b_{m,n}\}_{m,n\in\N}\in \M(\A)\}$. We shall examine the product
$$
c_{p,r,s,t}=U^{*p}\pi(a_{p,r})U^{r}  U^{*s}\pi(b_{s,t} )U^{t}
$$
Depending on the relationship between $r$ and $s$ we have two cases.
\begin{Item}{1)}
 If $s \leq r$ then
$$
c_{p,r,s,t}=U^{*p}\pi(a_{p,r})U^{r-s} (U^s U^{*s})\pi(b_{s,t}) U^t=
U^{*p}\pi(a_{p,r})U^{r-s} U^{*r-s} U^{r-s}\pi(b_{s,t}) U^t
$$
$$
=U^{*p}\pi(a_{p,r})U^{r-s} \pi(b_{s,t})(U^{*r-s} U^{r-s}) U^t=U^{*p}\pi(a_{p,r} \delta^{r-s}(b_{s,t})) U^{t+r-s}
$$
Putting $r-s=j$, $r=i$, $p=m$ and $t+r-s =n$ we get
$$
c_{p,r,s,t}=U^m \pi(a_{m,i} \delta^{j}(b_{i-j,n-j})) U^{n}
$$
and thus
$$
\sum_{s,r \in \N \atop s \leq r}  c_{m,r,s,n-r+s}=
\sum_{j=0}^\infty\sum_{i=j}^\infty  U^m \pi(a_{m,i} \delta^{j}(b_{i-j,n-j})) U^{n}=
U^m \pi(\big(a\cdot\sum_{j=0}^\infty  \Lambda^j(b)\big)_{m,n})  U^{n}
$$
\end{Item}
\begin{Item}{2)}
If $r <s$ then
$$
c_{p,r,s,t}=U^{*p} \pi(a_{p,r})(U^r  U^{*r}) U^{*s-r}\pi(b_{s,t}) U^l=
U^{*p}\pi( a_{p,r})U^{*s-r}U^{r-s}U^{*r-s}\pi(b_{s,t}) U^t
$$
$$
=U^{*p} (U^{*s-r}U^{s-r})\pi(a_{p,r})U^{*s-r}\pi(b_{s,t} )U^t=U^{*p+s-r} \pi(\delta^{s-r}(a_{p,r})b_{s,t}) U^t
$$
Putting $s-r=j$, $r=i$, $p +s-r=m$ and $t=n$ we get
$$
c_{p,r,s,t}=U^m\pi( \delta^{j}(a_{m-j,i-j})b_{i,n} ) U^{n}
$$
and thus
$$
\sum_{s,r \in \N \atop r< s}  c_{m-s+r,r,s,n}=
\sum_{j=1}^\infty\sum_{i=j}^\infty  U^m \pi(\delta^{j}(a_{m-j,i-j})b_{i,n} ) U^{n}=
U^m \pi(\big(\sum_{j=0}^\infty  \Lambda^j(a)\cdot b\big)_{m,n})  U^{n}
$$
\end{Item}
\noindent Using the formulas obtained in 1) and 2)  we have
$$
\Psi_{(\pi, U)}(a) \Psi_{(\pi, U)}( b)= \sum_{p,r,s,t\in \N} c_{p,r,s,t}=
\sum_{p,r,s,t\in \N \atop s\leq r}c_{p,r,s,t} + \sum_{p,r,s,t\in \N \atop r<s} c_{p,r,s,t}$$
$$
=\sum_{m,r,s,n\in \N \atop s\leq r,\,  n\leq r-s}c_{m,r,s,n-r+s} +
\sum_{m,r,s,n\in \N \atop r<s,\,m\leq s-r} c_{m-s+r,r,s,n}=
\sum_{m,n\in \N} U^{*m} (a\star b)_{m,n} U^n =\Psi_{(\pi, U)}(a\star b)
$$
and the proof is complete.
\end{Proof}

 We now examine the structure of $\M(\A)$.
    We shall say that  a matrix $\{a_{n}^{(m)}\}_{n\in\N, m\in \Z}$ in $\M(\A)$  is  $k$-\emph{diagonal},
    where $k$ is an integer, if it satisfies the condition
$$
a_{n}^{(m)}\neq 0  \Longrightarrow\,\, m=k.
$$
In other words $k$-diagonal matrix is the one of the form
\begin{center}\setlength{\unitlength}{1mm}
\begin{picture}(110,24)(-5,-10)

\put(-10,0){$
 \left(\begin{array}{c}\begin{xy}
\xymatrix@C=-1pt@R=3pt{
   &      \, \,   &    \qquad \,\,0 \,   \\
     &      \,     &       \\
      &     \, 0      &       \\
     &       &   \qquad
        }
  \end{xy}
  \end{array}\right)
 $  }
  \put(3.5,11){\scriptsize $k$}
  \qbezier[10](2, 10)(4,10)(6, 10)

 \qbezier[42](2, 10)(9,2.5)(16,-5)
  \qbezier[36](6, 10)(11,4.5)(16,-1)

   \put(29,0){if $k\geq 0$, or}
  \put(98,0){if $k< 0$.}
 \put(61,0){$
 \left(\begin{array}{c}\begin{xy}
\xymatrix@C=-1pt@R=4pt{
   &      \, \,   &    \,   \\
     &      \,     &  0      \\
      &     \,      &       \\
   0 \,\, &       &   \qquad
        }
  \end{xy}
  \end{array}\right)$}
  \put(58,3.6){\scriptsize $|k|$}
  \qbezier[10](64,2) (64,4) (64,5.5)
  \qbezier[42](64,5.5)(72.4,-1.5)(80.8,-8)
  \qbezier[36](64,2)(70,-3.5)(76,-8)
     \end{picture}
     \end{center}
  The  linear space consisting of all $k$-diagonal matrices will be denoted by $\mathcal{M}_{k}$.
  These will correspond to spectral subspaces.
\begin{prop}
\label{**}
The spaces $\M_k$ define  a $\Z$-gradated algebra structure on $\M(\A)$. Namely
$$
\M(\A)=\bigoplus_{k\in \Z} \M_k,$$
 and for every $k$, $l\in \Z$ we have the following relations
$$
\M_{k}^*=\M_{-k},\qquad\M_{k}\star\mathcal{M}_{l}\subset \M_{k+l}.
$$
In particular, $\M_{0}$ is a $^*$-algebra,  $\M_{k}\star \M_{-k}$
is a self-adjoint two sided  ideal in $\M_{0}$. Moreover
$$
 \M_k\star \M_k^*\star \M_k=\M_k, \qquad \M_k\star \M_k^*
 \star \M_l\star \M_l^*=\M_l\star \M_l^* \star \M_k\star \M_k^*.
$$
\end{prop}
\label{coef-alg}
\begin{Proof}
Relations  $\M_k^*=\M_{-k}$, $\M_k\star \M_l\subset \M_{k+l}$ and
$(\M_k \star \M_k^*)^*=(\M_k \star \M_k^*)$ can be checked by
means of an elementary matrix calculus. Using these relations we
get
$$
\M_0 \star (\M_k \star \M_k^*)= (\M_0 \star \M_k)\star \M_k^*\subset \M_k\star \M_k^* ,
$$
$$
(\M_k\star \M_k^*)\star \M_0 = \M_k\star (\M_k^*\star
\M_0)\subset  (\M_k \star \M_k^*),
$$
and thus $\M_k\star \M_k^*$ is an ideal in $\M_0$. Since the product of two ideals
 is equal to their intersection we have
 $\M_k\star \M_k^* \star \M_l\star \M_l^*= \M_l\star \M_l^* \star \M_k \star \M_k^*$.
 In order to see that $\M_k\star  \M_k^*\star \M_k= \M_k$
 one readily checks that  the identity $1_k$ in $\M_k \star \M_k^*$  satisfies  $1_k a=a$ for every  $a\in \M_k$.
\end{Proof}

Proposition \ref{**} indicates in particular that $\M_0$ may be  regarded as  a coefficient algebra for  $\M(\A)$.
This will be showed explicitly in Proposition \ref{nie dam rady umre}.

\begin{cor}\label{wniosek o spektralnych podprzestrzeniach}
  Let $(\pi,U,H)$ be  covariant representation of $(\A,\delta)$ and let
  $
B_k= \Psi_{(\pi, U)}(\M_k)
 $
 be the linear space consisting of the elements of the form
 $$
\sum_{n=0}^N U^{*n} \pi(a_{n}^{(k)}) U^{n+k}, \quad  \textrm{ if }\,\, k \geq 0, \quad \textrm{ or}
\quad \sum_{n=0}^N U^{*n+|k|} \pi(a_{n}^{(k)})U^{n}, \quad  \textrm{ if }\,\, k < 0.
$$
Then for every $k$ and $l\in \Z$ we have the following relations
$$
B_k^*=B_{-k},\qquad B_kB_l\subset B_{k+l}
$$
In particular, $B_0$ is a $C^*$-algebra,  $B_k B_k^*$ is a self-adjoint two sided  ideal in $B_0$ and
$$
B_k B_k^* B_k=B_k, \qquad B_k B_k^*  B_l B_l^*=B_lB_l^*  B_k B_k^*.
$$
\end{cor}
The importance of $B_0$ was observed in \cite{Leb-Odz} and is clarified  by the next proposition,
see \cite[Proposition 2.4]{Leb-Odz}.
\begin{prop}
Let $(\pi,U,H)$ be  covariant representation of $(\A,\delta)$ and adopt the notation
from Theorem \ref{dense subalgera structure prop} and Corollary \ref{wniosek o spektralnych podprzestrzeniach}.
Every element $a \in C_0^*(\pi(\A),U)$ can be presented in the form
$$
a= \sum_{k=1}^{\infty} U^{*k} a_{-k} + \sum_{k=0}^{\infty}  a_{k}U^{*k}
$$
where $a_{-k} \in B_0\pi(\delta^k(1))$,  $a_{k} \in \pi(\delta^k(1))  B_0$, $k\in \N$,
and only finite number of these coefficients  are non-zero.
\end{prop}
We shall now formulate a similar result concerning $\M(\A)$. For
each $k\in \Z$ we define a mapping $N_k:\M(\A)\to \M_0$, $k\in\Z$,
that carries a $k$-diagonal onto a $0$-diagonal and delete all the
remaining ones. Namely, for $a = \{ a_{n}^{(k)}  \}$ we set
$$
\left[N_k (a)\right]_{n}^{(m)}=
\begin{cases}
a_{n}^{(k)} & \textrm{ if }  m=0,\\
0 & \textrm{ otherwise },
\end{cases} \qquad \quad k \in \Z.
$$
One readily checks that   for $k\geq 0$ we have $N_k(\M_k)
=\M_0\star \delta^k(1)$, $N_{-k}(\M_{-k})= \delta^k(1) \star
\M_0$. Thus  we get that  the algebra  $\M_0$ consists of elements
that play the role of Fourier coefficients in $\M(\A)$.
\begin{prop}\label{nie dam rady umre}
Every element $a$ of  $\M(\A)$ is uniquely presented in the form
$$
a= \sum_{k=1}^{\infty} u^{*k}\star a_{-k} + \sum_{k=0}^{\infty}  a_{k}\star  u^{*k}
$$
where $u$ is given by (\ref{notation stupid}) and
 $a_{-k} \in \M_0\star \delta^k(1)$,  $a_{k}
\in \delta^k(1) \star \M_0$, $k\in \N$, and only finite number of
these coefficients  is non-zero. Namely, $a_k=N_k(a)$ for $k\in
\Z$.
\end{prop}

\section{Norm evaluation of elements in $C_0^*(\pi(\A),U)$}

In this section we gather a number of technical results concerning norm evaluation of elements in $C_0^*(\pi(\A),U)$.
We shall make use of these results in the subsequent sections.

The mappings $N_k:\M_k \to \M_0$ factors through $\Phi_{(\pi,U)}$ to mappings $\NN_k:B_k\to B_0$.
\begin{prop}\label{proposition for k-diagonals}
Let $(\pi,U,H)$ be  a covariant representation of $(\A,\delta)$ and let $k\in \Z$.  Then the norm $\|a\|$ of an element $a\in B_k$ corresponding to the matrix $\{a_{n}^{(m)}\}_{n\in \N, m\in \Z}$  in $\M_k$ is given by
$$
\lim_{n\to \infty} \max \left\{\max_{i=1,...,n}\Big \|(1-U^*U)\sum_{j=0}^{i}
\pi(\delta^{i-j}(a_{j}^{(k)}))\Big\|,\, \Big\|U^*U \pi(a_{n}^{(k)})\Big\| \right\}.
$$
In particular the mapping $\NN_k:B_k\to B_0$ given by
$$
\NN_k (\Phi_{(\pi,U)}(a))=\Phi_{(\pi,U)}(N_k (a)), \qquad a\in \M_k.
$$
is a well defined linear isometry  establishing the following isometric isomorphisms
$$
B_k\cong B_0\delta^k(1),\quad \textrm{ if }\,\, k\geq 0,\qquad B_{k}\cong \delta^{|k|}(1) B_0,\quad
\textrm{ if }\,\, k < 0.
$$
\end{prop}
\begin{Proof} Let us assume that $k \geq 0$.
Let $N$ be such that $a_{m}^{(k)}=0$ for $m>N$. Then, similarly  as it was done in \cite[Proposition 3.1]{kwa4},
one checks that defining
$$
a_i=(1-U^*U) \pi(\sum_{j=0}^{i} \delta^{i-j}(a_{j}^{(k)})), \quad  i=0,...,N,\qquad a_{N+1}= U^*U \pi(a_{N}^{(k)}),
$$
we have
 \begin{equation}\label{postac coefficientowa}
a = \big(a_0 +U^*a_1U +... + U^{*N}(a_N + a_{N+1})U^N\big) U^k
\end{equation}
where
\begin{equation}\label{coefficients warunki1}
a_i\in (1-U^*U)\pi(\delta^i(1)\A \delta^{i+k}(1)), \quad  i=0,...,N,\quad a_{N+1} \in U^*U\pi(\delta^N(1)\A
\delta^{N+k}(1)),
\end{equation}
Since $U^k$ is a partial isometry,  fomulae  \eqref{postac coefficientowa}  and \eqref{coefficients warunki1}
imply the following equalities
 $$
 \|a\|=\| (a_0 +U^*a_1U +... + U^{*N}(a_N+ a_{N+1}) U^N) U^kU^{*k}\|
 $$
 $$
 =\| a_0\pi(\delta^k(1)) +U^*a_1\pi(\delta^{k+1}(1))U +... + U^{*N}(a_N+a_{N+1})\pi(\delta^{k+N}(1)) U^N \|
 $$
 $$
 =\|a_0 +U^*a_1U +... + U^{*N}(a_N+ a_{N+1})U^N\|=\max_{i=0,...,N+1}\|a_i\|,
 $$
where the final equality follows from \eqref{coefficients warunki1} and
the fact that $U^{*j}U^{j}-U^{*j+1}U^{j+1}$, $j=0,..,N$, and $U^{*N+1}U^{N+1}$ are
pairwise orthogonal projections lying in $\pi(\A)'$, cf. \cite[Proposition 3.1]{kwa4}.
This proves the case when $k \geq 0$.\\
In the case of negative $k$ one may apply the part of proposition  proved above  to the adjoint $a^*$
of the element $a$ and thus the hypotheses follows.
\end{Proof}
We denote by
$$
d(a, K)=\inf_{b\in K} \|a-b\|
$$
the usual distance of an element $a$ from the set $K$. The definition of an ideal  associated
to a covariant representation (Definition \ref{kowariant  rep defn 2})  and a known  fact expressing
quotient norms in terms of projections, see for instance, \cite[Lemma 10.1.6]{Kadison},  gives us the following
\begin{cor}\label{corollary with distance}
If $(\pi,U,H)$ is a covariant representation of $(\A,\delta)$ associated with an ideal $J$
and $I$ denotes the kernel of $\delta$, then the norm of an element $a\in B_k$ corresponding to a matrix
$\{a_{n}^{(m)}\}_{n\in \N, m\in \Z}$  in $\M_k$ is given by
$$
\|a\|=\lim_{n\to \infty} \max \left\{\max_{i=1,...,n}\big\{ d\big(\sum_{j=0}^{i}
\delta^{i-j}(a_{j}^{(k)}),J\big)\big\},\, d(a_{n}^{(k)},I) \right\}.
$$
In particular, the spaces $B_k$, $k\in\Z$ do not depend on the choice of covariant representation associated
with  a fixed ideal $J$ orthogonal to $I$.
\end{cor}
We showed in Proposition \ref{proposition for k-diagonals} that for an arbitrary covariant representation $(\pi,U,H)$,
the mappings $N_k$ factor through to the mappings $\NN_k$ acting on spaces $B_k$, in general however,
 $N_k:\M(\A)\to B_0$  do not factor through $\Psi_{(\pi,U)}$ to the mappings acting on
 the algebra $C_0^*(\pi(\A),U)$. In fact, it is the case iff the covariant representation  $(\pi,U,H)$
 satisfies a  certain property we are just about to introduce.
\begin{defn} We shall say that a covariant representation $(\pi,U,H)$ of $(\A,\delta)$
possesses the {\em property} $(*)$ if for any $a\in C_0^*(\pi(\A),U)$ given by
a matrix $\{a_{mn}\}_{m,n\in\N}\in\M(\A)$  the inequality
$$
\|\sum_{m\in\N} U^{*m}\pi(a_{m,m})U^m \| \leq \|\sum_{m,n\in\N} U^{*m}\pi(a_{m,n})U^n \|,  \qquad\qquad (*)
$$
holds. In view of Corollary \ref{corollary with distance}  the above equality could be equivalently stated
in the form    $$
\lim_{n\to \infty} \max \left\{\max_{i=1,...,n}\big\{ d\big(\sum_{j=0}^{i}
\delta^{i-j}(a_{j,j}),J\big)\big\},\, d(a_{n,n},I) \right\} \leq \|a\|,\quad (*)
$$
where  $(\pi,U,H)$ is associated with an ideal $J$ and $I$ is the kernel of $\delta$.
\end{defn}
The next result, see \cite[Theorem 2.8]{Leb-Odz}, indicates that under the fulfillment of property (*)
elements of $B_0$ play the role of 'Fourier' coefficients in the algebra $C_0^*(\pi(\A),U)$.
\begin{thm}
Let $(\pi,U,H)$ be covariant representation possessing the  property $(*)$
then the mappings $\NN_k: C_0^*(\pi(\A),U)\to B_0$, $k\in \Z$,  given by formulae
\begin{equation}\label{rzuty na coefficienty}
 \NN_k( \Phi_{(\pi,U)}(a)) = \sum_{n\in\N} U^{*n}\pi(a_{n}^{(k)})U^{n},
\end{equation}
where $\{a_{n}^{(m)}\}_{n\in \N, m\in \Z}\in \M(\A)$, are well defined contractions and thus
they  extend uniquely to bounded operators on $C^*(\pi(\A),U)$. In particular,
every element $a \in C_0^*(\pi(\A),U)$ can  be uniquely presented in the form
$$
a= \sum_{k=1}^{\infty} U^{*k} a_{-k} + \sum_{k=0}^{\infty}  a_{k}U^{*k}
$$
where $a_{-k} \in B_0\pi(\delta^k(1))$,  $a_{k} \in \pi(\delta^k(1))  B_0$, $k\in \N$,
namely, $a_k=\NN_k(a)$, $k\in \Z$.
\end{thm}

Let  us also recall  \cite[Theorem 2.11]{Leb-Odz}.
\begin{thm}\label{3a.N}
If $(\pi,U,H)$ possess the  property $(*)$, then for any element $a$ in $C_0^*(\pi(\A),U)$  we have
\begin{equation}\label{be3.131}
\Vert a \Vert = \lim_{k\to\infty}
\sqrt[\leftroot{-2}\uproot{1}\scriptstyle 4k]{
\left\Vert \NN_0 \left[ (aa^*)^{2k}\right]\right\Vert }
\end{equation}
where $\NN_0$ is the mapping defined by \eqref{rzuty na coefficienty}.
\end{thm}

Using the above results  one sees that in the presence of property $(*)$ the  norm of an
element $a\in C_0^*(\pi(\A),U)$ may be calculated only in terms of the elements
of $\A$. Indeed, as $\NN_0 \left[ (aa^*)^{2k}\right]$ belongs to $B_0$ one can apply Corollary
\ref{corollary with distance} to calculate $\| \NN_0 \left[ (aa^*)^{2k}\right]\|$ in terms of
the matrix from $\M(\A)$ corresponding to $a$. However in practice, the  calculation of
the matrix corresponding to the element $(aa^*)^{2k}$ starting from $a$, see formula \eqref{star1},
seems to be an extremely difficult task.

\section{Crossed product}

Now we proceed to the description of the goal (the main object) of the article.

The set $\M(\A)$ with operations \eqref{add1},
\eqref{mulscal1}, \eqref{invol1}, \eqref{star1} is an algebra with involution. We define
a seminorm  on $\M(\A)$ that will depend on the choice of an orthogonal ideal. Let $J$ be
a fixed ideal in $\A$ having zero intersection with kernel of $\delta$.
Let
$$
    \|| a\||_J:=\sum_{k\in \Z} \lim_{n\to \infty} \max \left\{\max_{i=1,...,n}\big\{ d\big(\sum_{j=0}^{i}
    \delta^{i-j}(a_{j}^{(k)}),J\big)\big\},\, d(a_{n}^{(k)},I) \right\}
    $$
    where $a=\{a_{n}^{(k)}\}_{n\in \N,k\in \Z}\in \M(\A)$.
    \begin{prop}
The  function  $\|| \cdot \||_J$ defined above is a  seminorm on $\M(\A)$ which is  $^*$-invariant and submulitplicative.
\end{prop}
\begin{Proof}
Let $(\pi,U,H)$ be an arbitrary covariant representation of $(|A,\delta)$ associated with $J$.
Such representation does exist by Theorem \ref{existance theorem}.
Then in view of  Corollary \ref{corollary with distance}  for every $a\in \M_k$, $k\in \Z$, we have
$$
    \|| a\||= \| \Phi_{(\pi,U)}(a)\|
    $$
where $\Phi_{(\pi,U)}:\M(\A)\to L(H)$ is the $^*$-homomorphism defined in
Proposition \ref{dense subalgebra structure prop}.
Using the fact that every element $a\in \M(\A)$ can be presented in the
form $a= \sum_{k\in \Z} a^{(k)}$ where  $a^{(k)}\in \M_k$ one easily sees
that $\||| \cdot \|||$ is $^*$-invariant seminorm. To show that it is
submultiplicative take $a= \sum_{k\in \Z} a^{(k)}\in \M(\A)$ and  $b= \sum_{k\in \Z} b^{(k)}\in \M(\A)$
such that $a^{(k)}, b^{(k)}\in \M_k$.  Then
$$
\|| a \star b \|| = \||\sum_{k\in \Z} a^{(k)}\star \sum_{l\in \Z} b^{(l)}\||=\||\sum_{k\in \Z}
\sum_{l\in \Z} a^{(k)}\star b^{(l)}\||\leq \sum_{k,l\in \Z}  \||a^{(k)}\star b^{(l)}\||
$$
$$
= \sum_{k,l\in \Z}  \|\Phi_{(\pi,U)}(a^{(k)}\star b^{(l)})\| \leq  \sum_{k,l\in \Z}  \|\Phi_{(\pi,U)}(a^{(k)})\|
\cdot  \| \Phi_{(\pi,U)}(b^{(l)})\|
$$
$$
 = \sum_{k,l\in \Z}  \||a^{(k)} \|| \cdot \||b^{(l)}\||
=\sum_{k \in\Z} \||a^{(k)}\|| \sum_{l\in \Z} \||b^{(l)}\||  =\||a\|| \cdot\|| b\||.
$$
\end{Proof}

\begin{defn}
Let $(\A,\delta)$ be a $C^*$-dynamical system and $J$ an ideal in $\A$ having a zero intersection with the kernel
of $\delta$. Let
$$
C^*(\A,\delta,J)
$$
be the enveloping $C^*$-algebra  of the quotient $^*$-algebra $\M(\A)/ \|| \cdot \||_J$.
The $C^*$-algebra $C^*(\A,\delta,J)$ will be called a \emph{crossed product} of $\A$ by $\delta$ associated with $J$.
\end{defn}
 Regardless of $J$, composing the quotient map with natural embedding of $\A$ into $\M(\A)$
 one has an embedding of $\A$ into $C^*(\A,\delta,J)$. Moreover, denoting by $\hat{u}$
 an element of $C^*(\A,\delta,J)$ corresponding to $u\in \M(\A)$ (see (\ref{notation stupid})),
 one sees that $C^*(\A,\delta,J)$ is generated by $\A$ and $\hat{u}$.

 \section{Isomorphism theorem}

 Once a universal object (the crossed product) is defined it is reasonable to have its faithful representation.
 This section is devoted to the description of the properties of such representations.

\begin{thm}[\textbf{Isomorphism Theorem}]\label{isomorphiasm theorem}
\label{iso}
Let $J$ be an ideal in $\A$ having a zero intersection with the kernel $I$ of $\delta$ and
let $(\pi_i,U_i,H_i)$, $i=1,2$, be covariant representations of $(\A,\delta)$ associated with $J$
and possessing the property $(*)$.
Then the relations
$$ \Phi(\pi_1(a)):=\pi_2(a),\quad a\in \A,\qquad
\Phi(U_1):=U_2
$$
gives rise to the isomorphism between the $C^*$-algebras $C^*(\pi_1(\A),U_1)$
and $C^*(\pi_2(\A),U_2)$.
\end{thm}
\begin{Proof} Let $B_{0,i}$ be a $^*$-algebra consisting of elements of the form
$\sum_{n=0}^N U^{*n}_i\pi_i(a_n)U_i^n$, $i=1,2$.
In view of Corollary \ref{corollary with distance}, $\Phi$ extends to the isometric isomorphism
from $B_{0,1}$ onto $B_{0,2}$. Moreover, we have
$$
\Phi(U_1aU_1^*)=U_2(\Phi(a))U^*_2, \qquad a \in B_{0,1}.
$$
Hence the assumptions of \cite[Theorem 2.13]{Leb-Odz} are satisfied and the hypotheses follows.
\end{Proof}
\begin{cor}
If $(\pi,U,H)$ possess property $(*)$, then we have the action $\al$ of the group $S^1$ on $C^*(\pi(\A),U)$
by authomorphisms given by
$$
\al_z(\pi(a)):=\pi(a), \quad a\in \A, \qquad  \al_z(U) := z U, \qquad z\in S^1.
$$
Moreover the spaces $\overline{B}_k$ are the spectral subspaces corresponding to this action, that is we have
$$
\overline{B}_k=\{ a \in C^*(\A,U): \al_z(a)= z^k a\}.
$$
In particular, the $C^*$-algebra $\overline{B}_0$ is the fixed point algebra for $\al$.
\end{cor}
\begin{Proof} Let $(\pi,U,H)$ be associated with $J$ and let $z\in S^1$. It is clear that $(\pi,zU,H)$
is also a covariant representation of $(\A,\delta)$ associated with $J$ and $(\pi,zU,H)$ possess property $(*)$.
Hence by using Theorem \ref{isomorphiasm theorem}, $\al_z$  extends
to the isomorphism of $C^*(\pi(\A),U)=C^*(\pi(\A),zU)$. The remaining part of the statement is obvious.
\end{Proof}
The next  theorem is an immediate corollary of the previous statements and \cite[Theorem 2.15]{Leb-Odz}.
It is another manifestation of the fact that  the elements $\NN_k (a)$, $k\in\Z$, should be considered as
Fourier coefficients for $a \in C^*(\pi(\A),U)$.
\begin{thm}
\label{uniqueNk}
Let $(\pi,U,H)$ possess property $(*)$ and let
\begin{center}
$a \in C^*(\pi(\A),U)$
\end{center}
Then the following conditions are equivalent:

\smallskip
\quad\llap{$(i)$}\ \ $a=0;$

\smallskip
\quad\llap{$(ii)$}\ \ $\NN_k (a)=0$, \,$k\in\Z;$

\smallskip
\quad\llap{$(iii)$}\ \ $\NN_0 (a^*a)=0$.

\end{thm}

The results presented above give us a possibility to write out a
  criterion for the representation of the crossed product to be faithful.

\begin{thm} If  $(\pi,U,H)$, is covariant representations of $(\A,\delta)$ associated with $J$, then
relations
$$
(\pi\times U)(a)=\pi(a),\qquad (\pi\times U)(\hat{u})=U
$$
determines in unique way an epimorphism $\pi\times U: C^*(\A,\delta, J) \to C^*(\A,U)$. Moreover
$\pi\times U$ is an isomorphism iff $(\pi,U,H)$ possesses property $(*)$

\end{thm}

 \noindent \textsc{Institute of Mathematics,  University  of Bialystok},\\
\textsc{ ul. Akademicka 2, PL-15-267  Bialystok, Poland }\\
 \emph{e-mail:} \texttt{bartoszk@math.uwb.edu.pl}\\
 \emph{www:} \texttt{http://math.uwb.edu.pl/$\sim$zaf/kwasniewski}
 \bigskip

\noindent \textsc{Institute of Mathematics,  University  of Bialystok},\\
\textsc{ ul. Akademicka 2, PL-15-267  Bialystok, Poland }\\
 \emph{e-mail:} \texttt{lebedev@bsu.by}\\

\begin{thebibliography}{99}




 \bibitem[ABL05]{Ant-Bakht-Leb}
 A.B. Antonevich, V.I. Bakhtin, A.V. Lebedev, "Crossed product of $C^*$-algebra by an endomorphism,
  coefficient algebras and transfer operators",  arXiv:\penalty0
math.OA/0502415 v1 19 Feb 2005.







 \bibitem[BL05]{Bakht-Leb}
V.\,I.\ Bakhtin, A.\,V.\ Lebedev,  "When a  $C^*$-algebra is a
coefficient algebra for a given endomorphism", arXiv:\penalty0
math.OA/0502414 v1 19 Feb 2005.









\bibitem[CK80]{CK}
  J.\ Cuntz and W.\ Krieger,
  "A Class of C*-algebras and Topological Markov Chains",
  { Inventiones Math.\ {\bfseries 56}, (1980), p.~251--268}





\bibitem[Exel94]{exel1}
R.\ Exel,   "Circle actions on $C^*$-algebras, partial
automorphisms and generalized Pimsner-Voiculescu exact sequence",
J.\ Funct.\ Analysis {\bfseries 122} (1994), p.~361--401.

\bibitem[Exel03]{exel2} R. Exel:  "A new look at the crossed-product of a $C^*$-algebra by an endomorphism",
Ergodic Theory Dynam. Systems, Vol  23, (2003), pp. 1733-1750,


\bibitem[Kwa05]{kwa}
B.K. Kwa\'sniewski: "Covariance algebra of a partial dynamical system", CEJM, 2005, V.3, No 4,  pp. 718-765
arXiv:\penalty0 math.OA/0407352


\bibitem[Kwa07a]{kwa3} B.K. Kwa\'sniewski: "On transfer operators for $C^*$-dynamical systems",
preprint arXiv:\penalty0 math.OA/0703798 v1 27 Mar 2007

 \bibitem[Kwa07b]{kwa4} B.K. Kwa\'sniewski: "Extensions of $C^*$-dynamical systems
to systems with complete transfer operators",
preprint  arXiv:\penalty0 math.OA/0703800 v1 27 Mar 2007

\bibitem[KL03]{maxid} B.K. Kwa\'sniewski and A.V. Lebedev: "Maximal ideal space of a commutative coefficient algebra",
preprint arXiv:\penalty0 math.OA/0311416





\bibitem[LO04]{Leb-Odz}
A.\,V.\ Lebedev, A.\ Odzijewicz ''Extensions of
$C^*$-algebras by partial isometries'', Matemat.\ Sbornik, 2004.
V.~195, No 7, pp.~37--70 (Russian).

\bibitem[KR86]{Kadison} Kadison R.V., Ringrose J.R.
\emph{Fundamentals of the theory of operator algebras. Vol.2. Advanced theory}, Academic Press, 1986.

\bibitem[Mur96]{Murphy}
  G.\,J.\ Murphy,
  "Crossed products of C*-algebras by endomorphisms",
  Integral Equations Oper. Theory {\bfseries 24}, (1996), p.~298--319.






\bibitem[Pas80]{Paschke}
W.\,L.\ Paschke, "The crossed product of a $C^*$-algebra by an
endomorphism", Proceedings of the AMS, {\bfseries 80}, No 1, (1980),
p.~113--118.


\bibitem[Ped79]{Pedersen} G. K. Pedersen: \emph{$C^*$-algebras and their automorphism groups},
Academic Press, London, 1979.




\bibitem[Sta93]{Stacey}
 P.\,J.\ Stacey, "Crossed products of $C^*$-algebras by
 $^*$-endomorphisms", J.\ Austral.\ Math.\ Soc. Ser.~A {\bfseries 54},
 (1993), p.~204--212.



\end{thebibliography}
\end{document}